\documentclass[review]{elsarticle}

\usepackage{lineno,hyperref}
\usepackage{amsmath}

\journal{Journal of \LaTeX\ Templates}

\newtheorem{theorem}{Theorem}
\newtheorem{lemma}{Lemma}

\newtheorem{corollary}{Corollary}

\numberwithin{equation}{section}

\begin{document}

	\begin{frontmatter}
		
		\title{Reconstruction of a function defined on a sphere using the Funk transform}		
		\author{Rafik Aramyan}
		\address{}
		\fntext[myfootnote]{The work was supported by the Science Committee of RA, in the frames of the research project 21AG‐1A045}

		\ead{rafikaramyan@yahoo.com}
		
		\author[mymainaddress]{Institute of Mathematics of NAS RA, Yerevan, Armenia}
		
		
		\address[mymainaddress]{Russian Armenian University, Yerevan, Armenia}

		\begin{abstract}
		It is known that the Funk transform (the Funk-Radon transform)
		is invertible in the class of even (symmetric) continuous functions defined on the unit 2-sphere $\mathbf S^{2}$.
		In this article, for the reconstruction of 			$f\in C(\mathbf S^{2})$ (can be non-even), an additional condition (to reconstruct an odd function) is found, and the injectivity of the so-called two data Funk transform is considered.  An iterative  inversion formula of the transform is presented. Such inversions have theoretical significance in convexity theory, integral geometry and spherical tomography. Also, the Funk-Radon transform is used in Diffusion-weighted magnetic resonance imaging.
		\end{abstract}
		
		
	\end{frontmatter}
	
	{\tiny }
\section{Introduction}\label{sec1}

\noindent 
Let $\mathbf S^2$ be the unit 2-sphere in $\mathbf R^3$ (the 3-dimensional Euclidean space) with the center at the origin (the space of unit vectors). By
${\mathbf S}_\omega\subset{\mathbf S}^2$ we denote the great circle with pole at
$\omega\in{\mathbf S}^2$. For a continuous function $f \in
C(\mathbf S^2)$, by
$f_\omega(\varphi)$, $\varphi\in{\mathbf S}_\omega$ we denote the restriction of $f$
on the circle ${\mathbf S}_\omega$. We define

\begin{equation}\label{1.1}
	Ff(\omega)=\dfrac{1}{2\pi}\int_{{\mathbf  S}_\omega}f_\omega(\varphi) \,d\varphi, \quad
	\omega\in{\mathbf  S}^2\quad\text{and}\quad \varphi\in{\mathbf S}_\omega,
\end{equation}
where $d\varphi$ is the usual angular measure, 
normalized by $2\pi$. $Ff(\omega)$ is the integral (mean)
of f over the circle ${\mathbf S}_\omega$. 

\noindent	The classical Funk transform (the Funk-Radon transform on ${\mathbf S}^2$) takes function $f \in
C(\mathbf S^2)$ to its integrals over great circles
\begin{equation}\label{1.2}
	f \to Ff.
\end{equation}

\noindent	We call that $f \in
C(\mathbf S^2)$: 

is an even (symmetric) function if $f(\omega)=f(-\omega) $ for all $\omega\in\mathbf S^2$;

is an odd function if $f(\omega)=-f(-\omega)$ for all $\omega\in\mathbf S^2$.

\noindent The Funk transform annihilates all odd functions, and so it is natural to consider the case when $f$ is even. In that case, the Funk transform takes even (continuous) functions to even continuous functions defined on $\mathbf S^2$.

\noindent		It is well known that the Funk transform is invertible in the class of even continuous functions defined on the unit sphere $\mathbf S^{2}$ (as usual in  literature this case is considered). The inversion formula is known whenever $Ff(\omega)$ is even and sufficiently smooth (see \cite{Fu}, \cite{Hel}, \cite{GGG}, \cite{Nat}, \cite{Rub02}).

\noindent The inversion of the Funk transform is required in the problems of convexity theory, integral geometry, spherical tomography and  many other areas of mathematics. The Funk transform is related to the well-known cosine transform and intersection bodies in convex geometry (\cite{AgQu}, \cite{Ara03}, \cite{B}-\cite{CH}, \cite{Gin}, \cite{Gar}, \cite{Hel}-\cite{Nat}, \cite{Pal}-\cite{Rub13}, \cite{WS}). Also, the Funk-Radon transform is used in Diffusion-weighted magnetic resonance imaging  (\cite{T}).

\noindent  A generalization of the Funk transform to an arbitrary dimension is due to Helgason \cite{Hel} (see also \cite{Rub13}). Several other  generalized Funk transforms were considered in  (\cite{Ag}-\cite{AgRu2}, \cite{Ara01}-\cite{Ara10}, \cite{Nat}, \cite{Gin}, \cite{Q}-\cite{Rub13}). The shifted Funk–Radon transform and non-geodesic Funk transform on the unit sphere were considered in (\cite{Ag}-\cite{AgRu2}, \cite{Q}-\cite{Q2}, \cite{Sal}-\cite{Sal2}).

\noindent It is clear that the Funk transform is non-invertible in the space of all continuous functions in $\mathbf S^{2}$ because it annihilates all odd functions.

\noindent \textbf{The problem arises}: finding an additional condition (to reconstruct an odd function) that allows us to recover an unknown function, using the Funk transform,  on the space of all continuous functions defined on $\mathbf S^{2}$. 	

\noindent In this article, for the reconstruction of
$f\in C(\mathbf S^2)$ using the Funk transform, the injectivity of the so-called two data Funk transform is considered.  Also, an inversion formula of the transform is presented.  Here, we apply the consistency method which was already used to invert the spherical Radon transform (SRT) in $\mathbf R^2$ (\cite{Ara19}) and $\mathbf R^3$ (\cite{Ara231}) (see also \cite{Ara10}). A similar problem was already solved in $\mathbf R^2$ to recover an odd function with respect to a line using SRT when the detectors are placed on the line (\cite{Ara232}). 

\noindent Throughout the article we use the
standard spherical coordinates $\omega=(\nu,\tau)$ for $\omega\in\mathbf S^2$ determined by the choice of the North Pole $\emph{N}\in\mathbf S^2$ and the reference
point $\tau=0$ on the equator ${\mathbf S}_{\emph{N}}$, $\nu\in[0,\pi]$  is  measured from the pole. Also, for $\omega=(\nu,\tau)$, where $\nu\in[0,\pi/2]$, we choose the direction  $\emph{N}_\omega=(\pi/2-\nu,\tau+\pi)$ as the reference point on $\mathbf S_\omega$ (the direction of the projection of the North Pole on $\mathbf S_\omega$) and let $\varphi\in[-\pi,\pi]$ determine the angular coordinate on   $\mathbf S_\omega$, the clockwise direction is positive. 

\noindent Given a continuous function 
$f\in {\textit{C}}(\mathbf S^{2})$, we define the following two weighted Funk transforms
\begin{equation}\label{1.3}
	Cf(\omega)=\dfrac{1}{2\pi}\int_{{\mathbf  S}_\omega}\cos\varphi \,\,f_\omega(\varphi) \,d\varphi=\dfrac{1}{2\pi}\int_{-\pi}^\pi\cos\varphi \,\,f_\omega(\varphi) \,d\varphi, 
\end{equation}
and
\begin{equation}\label{1.3.1}
	Sf(\omega)=\dfrac{1}{2\pi}\int_{{\mathbf  S}_\omega}\sin\varphi \,\,f_\omega(\varphi) \,d\varphi=\dfrac{1}{2\pi}\int_{-\pi}^\pi\sin\varphi \,\,f_\omega(\varphi) \,d\varphi, 
\end{equation}
here $f_\omega(\varphi)$, $\varphi\in{\mathbf S}_\omega$ is the restriction of $f$
on the circle ${\mathbf S}_\omega$ and we integrate with
respect to the normalized angular measure on ${\mathbf S}_\omega$. Note that the definition of $Cf(\omega)=C_\emph{N}f(\omega)$ (respectively $Sf(\omega)=S_\emph{N}f(\omega)$) depends on $\emph{N}$. In the case where we consider the transforms with respect to another direction, say $\Omega\in\mathbf S^{2}$ (not $\emph{N}$), we will write $C_\Omega f(\omega)$.

\noindent We denote by ${ \cal C}^{1}(\mathbf S^2)$ the space of all functions with continuous partial derivatives of first order.

\noindent \emph{The main results that we will prove in the following sections}.

\begin{theorem}\label{1} Let $f\in{ \cal C}^{1}(\mathbf S^2)$ (can be non even) and $\Omega\in\mathbf S^2$. The value of $f(\Omega)$ can be recovered using the two data Funk transform
	\begin{equation}\label{1.4}
		f\to (Ff(\omega),C_\Omega f(\omega))\,\,\, \text{for}\,\,\,\omega\in\mathbf S^2\end{equation}
	($\Omega\in\mathbf S^2$ is taken for the north Pole). 
\end{theorem}

\begin{figure}
	\center
	\includegraphics{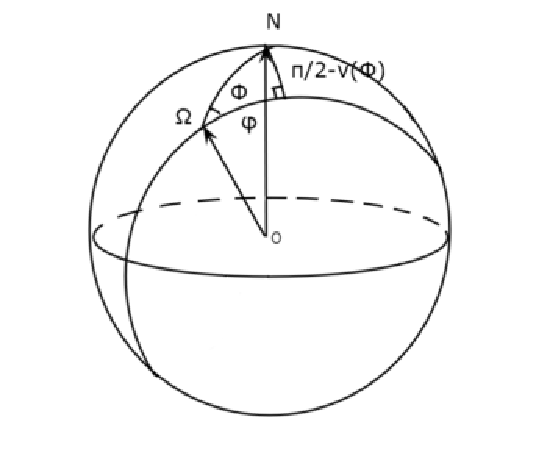}
	\caption{An illustration to (\ref{4.1}).}
\end{figure}

\noindent For $f\in{C}^{1}(\mathbf S^2)$	the Fourier expansion of $f_{\omega}(\varphi)$ we write as follows: 
\begin{equation}\label{1.5}
	f_{\omega}(\varphi)=Ff(\omega)+\sum_{n=1}^{\infty} \left[a_{n}(\omega) \,\cos(n\varphi)+b_{n}(\omega) \,\sin(n\varphi)\right],\end{equation}
here $Ff(\omega)=a_{0}(\omega)$,  $a_{n}(\omega)$ and $b_{n}(\omega)$ are the Fourier coefficients of the restriction  $f_\omega(\varphi)$ of $f$ on ${\mathbf S}_\omega$ and depend on the North Pole $\emph{N}$.

\noindent From (\ref{1.5}) it follows that 
\begin{equation}Ff(\omega)=a_{0}(\omega),\,\,\,\,\,\,\,\,Cf(\omega)=\dfrac{1}{2}a_{1}(\omega)\label{1.6}.\end{equation}

\noindent Thus, Theorem 1 states that to recover $f(\Omega)$ for $\Omega\in\mathbf S^2$ it is enough to have the Fourier first two coefficients $Ff(\omega)=a_{0}(\omega)$ and $a_{1}(\omega)$ of the restrictions of $f$ onto ${\mathbf S}_\omega$ ($\Omega\in\mathbf S^2$ is taken for the north Pole). We will prove Theorem 1 by presenting the reconstruction formula (see \ref{1.10}).

\noindent	It is well known that any $f$ defined on ${\mathbf S}^2$ can be represented as the sum of the even and odd functions as follows:

\begin{equation}f(\omega)=\dfrac{f(\omega)+f(-\omega)}{2}+\dfrac{f(\omega)-f(-\omega)}{2}=f^+(\omega)+f^-(\omega)\label{1.7}\end{equation}
where $f^+(\omega)$ is an even function and $f^-(\omega)$ is an odd function.

\noindent	Let $a^+_{n}((\nu,\tau))$ and $b^+_{n}((\nu,\tau))$ (respectively $a^-_{n}((\nu,\tau))$ and $b^-_{n}((\nu,\tau))$) are the Fourier coefficients of the restriction  $f^+_\omega(\varphi)$ (respectively for  $f^-_\omega(\varphi)$).

\noindent	We have the following identities $a^+_{1}((\nu,\tau))\equiv0$, $a^-_{0}((\nu,\tau))\equiv0$ and 
\begin{equation}a^+_{0}((\nu,\tau))=a_{0}((\nu,\tau)),\,\,\,\,\,\,\,\,\,\,
	a^-_{1}((\nu,\tau))=a_{1}((\nu,\tau))	
	\label{1.8}\end{equation}
hence, knowing $Ff(\omega)=a_{0}(\omega)$ of  an unknown function $f$, one can recover the even part $f^+(\omega)$ and  respectively, knowing $C_\Omega f(\omega)=a_{1}(\omega)/2$ for every $\Omega\in{\mathbf S}^2$, one can recover the odd part $f^-(\omega)$ of $f$.

\noindent \textbf{Problem 2}. An inversion formula.

\noindent We define two new sequences of trigonometric 
polynomials 	\begin{equation}\label{1.9}
	P_n^1(u)= \sum_{k=1}^{n}\sum_{m=1}^{k-1}c_m(2k-1){\sin^{2m}u},\,\,\,\text{and}\,\,\,P_n^0(u)=\sum_{k=1}^{n}\sum_{m=1}^{k}c_m(2k){\sin^{2m-1}u}\end{equation}
where, $c_m(2k-1)$ and $c_m(2k)$ are coefficients that one can calculate recursively (see (\ref{6.2})-(\ref{6.4}),  (\ref{6.5})-\ref{6.7}).
Polynomials are standard in the sense that their construction does not depend on $f$.

\noindent Here we present a reconstruction formula for $f(\emph{N})$. To recover $f(\Omega)$ for $\Omega\in\mathbf S^2$, one can choose $\Omega\in{\mathbf S}^2$ as the North Pole. Let 

\noindent $(Ff(\omega),C_\emph{N} f(\omega))=(Ff(\nu,\varphi),C f(\nu,\varphi))$ be the two data Funk transform  of $f$ defined on ${\mathbf S}^2$ (here we consider the spherical coordinates with respect to $\emph{N}$). For $\nu\in[0,\pi/2]$ we define the following functions
\begin{equation}\label{1.9.1}
	\overline{Ff(\nu)}=\int_0^{2\pi}F f(\nu,\tau)\,d\tau.\end{equation} and
\begin{equation}\label{1.9.2}
	\overline{Cf(\nu)}=\int_0^{2\pi}C f(\nu,\tau)\,d\tau.\end{equation}
and call them the averages of $Ff$ and $Cf$.

\begin{theorem}\label{thm2} Let $(Ff,Cf)$ be the two data Funk-Radon transform of a function $f\in{\textit{C}}(\mathbf S^2)$, we have 
	\begin{align}
		2\pi\,f(\emph{N})=\overline{Ff(\pi/2)}+\,\,\,\,\,\,\nonumber\\
		\lim_{n\to \infty}\left(2n(\overline{Cf(\pi/2)}+2\overline{Ff(\pi/2)})+\int_0^{\pi/2}\left( 2P_n^1(u)\overline{Cf(u)}+ P_n^0(u)\overline{Ff(u)}\right)\cos u\,du\right)
		\label{1.10}		\end{align}
	where $\overline{Ff(u)}$ and $\overline{Cf(u)}$ are the averages of the two data Funk transform of $f$ with respect to $\emph{N}$.		
\end{theorem}

\noindent Note that, the expression under the limit of (\ref{1.10}) represents a partial sum of a convergent series (see (\ref{5.2}) below).

\noindent Thus, to reconstruct $f\in{ C}(\mathbf S^2)$ we have to know the Fourier first two coefficients $a_{0}(\omega)=Ff(\omega)$ and $a_{1}(\omega)=2C_\Omega f(\omega)$ of $f$ for all $\Omega\in\mathbf S^2$ ($\Omega\in\mathbf S^2$ is taken for the north Pole).

\noindent	Knowing the two weighted Funk transforms $Cf(\omega)$ and $Sf(\omega)$ for a fixed direction (for the North Pole, for example), one can calculate $C_\Omega f(\omega)$ for all the directions $\Omega\in\mathbf S^2$ and $\omega\in\mathbf S^2$ (see Fig.2).

\begin{theorem}\label{thm2} Let $Cf(\omega)$ and $Sf(\omega)$ be the two weighted Funk transforms of a function $f\in{\textit{C}}^{1}(\mathbf S^2)$ with respect to the North Pole. For $\Omega\in\mathbf S^2$, we have 
	\begin{equation}
		C_\Omega f(\omega)=\cos\alpha\, Cf(\omega)+\sin\alpha \,Sf(\omega) 
		\label{1.11}		\end{equation}
	where $\alpha$ is the angular coordinate of the direction of the projection of $\Omega$  on ${\mathbf S}_\omega$ with respect to the direction of the projection of $\emph{N}$  on ${\mathbf S}_\omega$.		
\end{theorem}

\noindent  From Theorem 3 it follows that 

\begin{corollary} Let $Cf(\omega)$ and $Sf(\omega)$ are  the two weighted Funk transforms of a function $f\in{\textit{C}}^{1}(\mathbf S^2)$ with respect to the North Pole. The transform 	
	\begin{equation}\label{1.12}
		f \to 	(Ff(\omega),Cf(\omega), Sf(\omega))\,\,\,\,\emph{for}\,\,\,\omega\in\mathbf S^2\end{equation}
	is invertible in ${ \cal C}(\mathbf S^2)$. 
\end{corollary}

\section{Inversion formula in terms of the Fourier coefficients}\label{3}

\noindent Let $\emph{N}\in\mathbf S^2$ be the North Pole. For the angular coordinate of the reference direction $\emph{N}_\omega$ on ${\mathbf S}_\omega$  (the projection direction of $N$ to ${\mathbf S}_\omega$)  we have 
$\varphi=0$. 
Hence, using (\ref{1.5}), we get
\begin{equation}\label{3.1}
	f(\emph{N}_\omega)=f_{\omega}(0)=Ff(\omega)+\sum_{n=1}^{\infty}a_{n}(\omega).\end{equation}
where for $n\geq1$
\begin{equation}\label{3.1}
	a_{n}(\omega)=a_{n}(\nu,\tau)=\dfrac{1}{\pi}\int_{-\pi}^\pi f_\omega(\varphi)\,\cos {n\varphi} \,d\varphi.\end{equation}

\section{{The derivation of the Fourier coefficients}}\label{sec3}

\noindent \noindent We will derive the functions $a_n(\omega)$ and $b_n(\omega)$ from the consistency condition (see \cite{Ara10}, \cite{Ara19}, \cite{Ara231}). For $\Omega\in {\mathbf S}^2$, we consider the bundle of great circles containing $\Omega$. 	
Any great circle from the bundle is determined by the angular coordinate $\phi\in\mathbf S_\Omega$. Let  $(\nu(\phi),\tau(\phi))$ be the normal of the great circle from the bundle determined by $\phi\in\mathbf S_\Omega$ and $\varphi(\phi)$ determines $\Omega\in S_{(\nu(\phi),\tau(\phi))}$ (we will write $\Omega=(\nu(\phi),\tau(\phi),\varphi(\phi))$). 

\noindent	 The consistency method is based on the following statement. The restrictions  

\noindent	$f_{\omega}(\varphi)$ are consistent. We mean the following for $\Omega\in {\mathbf S}^2$, we have

\begin{equation}\label{4.0}
	f(\Omega)=f_{\omega}(\varphi) \,\,\,\,\emph{for all}\,\,\,\,\, \Omega=(\nu(\phi),\tau(\phi),\varphi(\phi))\end{equation}
hence, (\ref{4.0}) does not depend on $\phi$ (the consistency condition).

\noindent Let $\Omega=(\nu_\Omega,\tau_\Omega)\in {\mathbf S}^2$  ($\nu_\Omega\in (0,\pi/2)$), using spherical trigonometry we get (see Fig.1)
\begin{equation}\label{4.1}\begin{cases}	\cos\nu_\Omega=\cos\varphi(\phi)\,\sin\nu(\phi)\\
		-\sin\varphi(\phi)=\sin\tau(\phi)\,\sin\nu_\Omega\\
		\cos\nu(\phi)=\sin\phi\,\sin\nu_\Omega.\end{cases}
\end{equation}

\noindent Differentiating (\ref{4.1}) with respect to $\phi$, we get
\begin{equation}\label{4.2}\begin{cases}	\sin\varphi\,\varphi'_\phi\,\sin\nu=\cos\varphi\,\cos\nu\,\nu'_\phi\\
		-\cos\varphi\,\varphi'_\phi=\sin\nu_\Omega\,\cos\tau\,\tau'_\phi\\
		-\sin\nu\,\nu'_\phi=\cos\phi\,\sin\nu_\Omega.\end{cases}
\end{equation}
\noindent From (\ref{4.2}) we get ($\nu(\phi)\ne0$ since $\nu_\Omega\in (0,\pi/2)$)

\begin{equation}\label{4.4}
	\nu'_\phi=-\sin\varphi, \quad \tau'_\phi=-\frac{\cos\varphi}{\sin\nu},\quad \varphi'_\phi=-\frac{\cos\varphi\cos\nu}{\sin\nu}.
\end{equation}

\noindent Differentiating  (\ref{4.0}) we obtain the following Lemma.
\begin{lemma}\label{1} The following equation for the restrictions $f_{(\nu,\tau)}(\varphi)=f(\nu,\tau,\varphi)$ is valid
	\begin{equation}\label{4.5}
		(f(\Omega))'_\phi=
		-f'_\nu\,\sin\varphi-f'_\varphi\frac{\cos\varphi\cos\nu}{\sin\nu}-
		f'_\tau\frac{\cos\varphi}{\sin\nu}=0.
	\end{equation}
	
\end{lemma}

\noindent In order to get differential equations for   $a_n((\nu,\tau))$ and $b_n((\nu,\tau))$, we
multiply (\ref{4.5}) by $\sin{n\varphi}$, for $n\geq0$, and integrate with respect to $\varphi$ over $[-\pi,\pi]$. We obtain

\begin{equation}\label{4.6}
	\int_{-\pi}^\pi f'_\nu\,\sin\varphi\sin{n\varphi}\,d\varphi+\int_{-\pi}^\pi f'_\varphi\frac{\cos\varphi\sin{n\varphi}\cos\nu}{\sin\nu}\,d\varphi+\int_{-\pi}^\pi
	f'_\tau\frac{\cos\varphi\sin{n\varphi}}{\sin\nu}\,d\varphi=0.
\end{equation}

\noindent	Considering (\ref{4.6}) term by term, using Integration by parts, for $n\geq0$ we get. For $n=1$
\begin{equation}
	(a_2(\nu,\tau))'_\nu+2a_2(\nu,\tau)\cot\nu-2(a_0(\nu,\tau))'_\nu-\dfrac{b_2(\nu,\tau))'_\tau}{\sin\nu}=0 \label{4.7}\end{equation}
and for $n>1$
\begin{eqnarray}
	(a_{n+1}(\nu,\tau))'_\nu-(a_{n-1}(\nu,\tau))'_\nu+(n+1)a_{n+1}(\nu,\tau)\cot\nu+\nonumber\\(n-1)a_{n-1}(\nu,\tau)\cot\nu-\dfrac{(b_{n+1}(\nu,\tau))'_\tau}{\sin\nu}-\dfrac{(b_{n-1}(\nu,\tau))'_\tau}{\sin\nu}=0	\,\,.\label{4.8}\end{eqnarray}

\section{{Averaging the inversion formula}}\label{sec4}

\noindent  Now, consider the bundle of great circles containing $\emph{N}$. Let  $(\pi/2,\tau)$ be the normal of the great circle from the bundle determined by $\tau\in\mathbf S_\emph{N}$ and $\varphi=0$ determines $\emph{N}\in S_{(\pi/2,\tau)}$. 

\noindent Using (\ref{4.1}), we get
\begin{equation}\label{5.1}
	f(\emph{N})=f_{(\pi/2,\tau)}(0)=Ff((\pi/2,\tau))+\sum_{n=1}^{\infty}a_{n}((\pi/2,\tau)).\end{equation}

\noindent We integrate both sides of \eqref{5.1} with respect to the arc
measure $d\tau$ over $[0,2\pi)$. We have (see \eqref{1.10}, \eqref{1.11} )

\begin{equation}\label{5.2}
	2\pi f(\emph{N})=\int_0^{2\pi}Ff((\pi/2,\tau))\,d\tau+\sum_{n=1}^{\infty}\int_0^{2\pi}a_{n}((\pi/2,\tau))\,d\tau=\overline{Ff(\pi/2)}+\sum_{n=1}^{\infty}\overline{a_n(\pi/2)}\end{equation}
where 	$\overline{a_n(\pi/2)}=\int_0^{2\pi}a_{n}((\pi/2,\tau))\,d\tau$.

\noindent 	Now the problem is to calculate $\overline{a_n(\pi/2)}$ for $n\geq1$.

\noindent Integrating both sides of \eqref{4.7} and \eqref{4.8} and taking into account
that
$$
\int_0^{2\pi}(b_{n}(\nu,\tau))'_\tau\,d\tau=0
$$
we get: for $n=1$
\begin{equation}
	(\int_0^{2\pi}a_2(\nu,\tau)\,d\tau)'_\nu+2\int_0^{2\pi}a_2(\nu,\tau)d\tau=2(\int_0^{2\pi}a_0(\nu,\tau)d\tau)'_\nu \label{5.7}\end{equation}
and for $n>1$
\begin{eqnarray}
	(\int_0^{2\pi}a_{n+1}(\nu,\tau)d\tau)'_\nu-(\int_0^{2\pi}a_{n-1}(\nu,\tau)\,d\tau)'_\nu+(n+1)\int_0^{2\pi}a_{n+1}(\nu,\tau)\,d\tau\,\cot\nu\nonumber\\+(n-1)\int_0^{2\pi}a_{n-1}(\nu,\tau)\,d\tau\,\cot\nu=0	\,\,.\label{5.8}\end{eqnarray}

\noindent Thus we have the following systems of differential equations

\begin{eqnarray}\label{5.10}\,\,\,\,\,\begin{cases}
		(\overline{a_2(\nu)})'_\nu+2	\overline{a_2(\nu)}\cot\nu=2(	\overline{a_0(\nu)})'_\nu\\
		(\overline{a_3(\nu)})'_\nu-(\overline{a_1(\nu)})'_\nu+3\overline{a_3(\nu)}\cot\nu+\overline{a_1(\nu)}\cot\nu=0\\
		(\overline{a_4(\nu)})'_\nu-(\overline{a_2(\nu)})'_\nu+4\overline{a_4(\nu)}\cot\nu+2\overline{a_2(\nu)}\cot\nu=0\\
		(\overline{a_{5}(\nu)})'_\nu-(\overline{a_{3}(\nu)})'_\nu+5\overline{a_{5}(\nu)}\cot\nu+3\overline{a_{3}(\nu)}\cot\nu=0\\
		\quad\cdots \\
		(\overline{a_{n+1}(\nu)})'_\nu-(\overline{a_{n-1}(\nu)})'_\nu+(n+1)\overline{a_{n+1}(\nu)}\cot\nu+(n-1)\overline{a_{n-1}(\nu)}\cot\nu=0\\\,\, \texttt{for}\,\,\,n\geq2,
\end{cases}\end{eqnarray}
where 

\begin{equation}\label{5.11}
	\overline{a_n(\nu)}=\int_0^{2\pi}a_{n}(\nu,\tau)\,d\tau.\end{equation}

\noindent For $\nu=0$ and $n\geq1$ we have the following boundary conditions 
\begin{equation}\label{5.12}
	\overline{a_n(0)}=\int_0^{2\pi}a_{n}(0,\tau)\,d\tau=\int_0^{2\pi}\dfrac{1}{\pi}\int_{S_{(0,\tau)}} f_{(0,\tau)}(\varphi)\,\cos {n(\tau-\varphi)}\,d\varphi \,d\tau=0.\end{equation}

\noindent
Finally to find $\overline{a_n(\nu)}$ for any $n>1$, we apply the following algorithm.  

\noindent	First, taking into account the boundary condition, one can find $\overline{a_2(\nu)}$ from the first equation of (\ref{5.10}). Substituting it into the third equation of (\ref{5.10}), one can calculate $\overline{a_4(\nu)}$ and so on, step by step one can calculate $\overline{a_n(\nu)}$ for all even $n$.

\noindent Also, separately one can calculate $\overline{a_3(\nu)}$ from the second equation of (\ref{5.10}) (note that $\overline{a_1(\nu)}$ is known). Substituting $\overline{a_3(\nu)}$ into the fourth equation of (\ref{5.10}), we obtain $\overline{a_5(\nu)}$ and so on, step by step one can calculate all coefficients $a_{n}$ for all odd $n>1$.
Now we are going to implement this algorithm.

\section{{A representation for the Fourier coefficients}}\label{sec5}

\noindent For $\overline{a_{2k+2}(\nu)}$,

\noindent ($\nu\in(0,\pi/2]$) we get (the  solution of (\ref{5.10}) with the boundary conditions (\ref{5.12}))

\begin{eqnarray}
	\overline{a_{n+1}(\nu)}=
	\frac{1}{\sin^{n+1}\nu} \int_{0}^{\nu}\sin^{n+1}u\,(\overline{a_{n-1}(u)})'_u-(n-1)\sin^{n}u\,\cos u\,\overline{a_{n-1}(u)}\, du=\nonumber\\\overline{a_{n-1}(u)}+\frac{1}{\sin^{n+1}\nu} \int_{0}^{\nu}-2n\sin^{n}u\,\cos u\,\overline{a_{n-1}(u)}\,du.\label{6.1}
\end{eqnarray}

\begin{theorem}\label{2} The following representations are valid. For $k\geq1$ we have		
	\begin{equation}\label{6.2}
		\overline{a_{2k}(\nu)}=2\overline{a_{0}(\nu)}+\frac{1}{\sin\nu}\int_0^\nu (\sum_{m=1}^{k}c_m(2k)\,(\dfrac{\sin u}{\sin\nu})^{2m-1})\cos u\,\overline{a_{0}(u)}\, du
	\end{equation}
	where $c_m(2k)$, $1\leq m\leq k$ are coefficients that one can calculate recursively.
\end{theorem}
\textbf{Proof.}  It follows from (\ref{5.10}) that  (\ref{6.2}) is true for $k=1$. Indeed
\begin{equation}\label{6.3}
	\overline{a_{2}(\nu)}=2\overline{a_{0}(u)}+\frac{1}{\sin\nu}\int_0^\nu -4\,(\dfrac{\sin u}{\sin\nu})\cos u\,\overline{a_{0}(u)} du
\end{equation}

\noindent	Suppose (\ref{6.2}) is true for some $n = k$. 	
Substituting (\ref{6.2}) into (\ref{6.1}), written for $n=2k+1$. We change the order of the summation, then the order of the integration, and after grouping them, we get

\begin{eqnarray}
	\overline{a_{2k+2}(\nu)}=
	2\overline{a_{0}(\nu)}+\frac{1}{\sin\nu}\int_0^\nu (\sum_{m=1}^{k}c_m(2k)\,(\dfrac{\sin u}{\sin\nu})^{2m-1})\cos u\,\overline{a_{0}(u)}\, du+\label{6.4}\\\frac{1}{\sin^{2k+2}\nu} \int_{0}^{\nu}(-4k-2)\sin^{2k+1}u\,\cos u\,(2\overline{a_{0}(u)})\,du+\nonumber\\\frac{(-4k-2)}{\sin^{2k+2}\nu} \int_{0}^{\nu}\sin^{2k}u\,\cos u\,\int_0^u (\sum_{m=1}^{k}c_m(2k)\,(\dfrac{\sin v}{\sin u})^{2m-1})\cos v\,\overline{a_{0}(v)}\, dv\,du=\nonumber\\2\overline{a_{0}(\nu)}+\frac{1}{\sin\nu}\int_0^\nu (\sum_{m=1}^{k}c_m(2k)(\dfrac{\sin u}{\sin\nu})^{2m-1}-(4k+2)(\dfrac{\sin u}{\sin\nu})^{2k+1})\cos u\,\overline{a_{0}(u)} du\nonumber\\+\frac{1}{\sin\nu}\int_0^\nu(\sum_{m=1}^{k}\dfrac{(-4k-2)c_m(2k)}{2k-2m+2}[(\dfrac{\sin u}{\sin\nu})^{2m-1}-(\dfrac{\sin u}{\sin\nu})^{2k+1}]\cos u\,\overline{a_{0}(u)} du.\nonumber
\end{eqnarray}

\noindent It follows from (\ref{6.4}) that (\ref{6.2}) is true for every $k$. Theorem 4 is proved.

\noindent For $\overline{a_{2k+1}(\nu)}$, (($\nu\in(0,\pi/2]$)) we get (the  solution of (\ref{5.10}) with the boundary conditions (\ref{5.12}))

\begin{eqnarray}
	\overline{a_{2k+1}(\nu)}=\label{6.4.1}\\
	\frac{1}{\sin^{2k+1}\nu} \int_{0}^{\nu}\sin^{2k+1}u\,(\overline{a_{2k-1}(u)})'_u-(2k-1)\sin^{2k}u\,\cos u\,\overline{a_{2k-1}(u)}\, du.\nonumber
\end{eqnarray}

\begin{theorem}\label{2} The following representations are valid. For $k>1$ we have		
	\begin{equation}\label{6.5}
		\overline{a_{2k-1}(\nu)}=	\overline{a_{1}(\nu)}+\frac{1}{\sin\nu}\int_0^\nu (\sum_{m=1}^{k-1}c_m(2k-1)\,(\dfrac{\sin u}{\sin\nu})^{2m})\cos u\,\overline{a_{1}(u)}\, du
	\end{equation}
	where $c_m(2k-1)$, $1\leq m\leq k$ are coefficients that one can calculate recursively.
\end{theorem}

\textbf{Proof.}  It follows from (\ref{5.10}) that  (\ref{6.5}) is true for $k=1$. Indeed
\begin{equation}\label{6.6}
	\overline{a_{3}(\nu)}=
	\overline{a_{1}(\nu)}+\frac{1}{\sin\nu}\int_0^\nu -4\,(\dfrac{\sin u}{\sin\nu})^2\,\cos u\,\overline{a_{1}(u)}\, du
\end{equation}

Substituting (\ref{6.5}) into (\ref{6.4.1}), written for $n=2k$. We change the order of the summation, then the order of the integration, and after grouping them, we get

\begin{eqnarray}
	\overline{a_{2k+1}(\nu)}=
	\overline{a_{1}(\nu)}+\frac{1}{\sin\nu}\int_0^\nu (\sum_{m=1}^{k-1}c_m(2k-1)\,(\dfrac{\sin u}{\sin\nu})^{2m})\cos u\,\overline{a_{1}(u)}\, du\label{6.7}\\+\frac{1}{\sin^{2k+1}\nu} \int_{0}^{\nu}-4k\sin^{2k}u\,\cos u\,\overline{a_{1}(\nu)}\,du+\frac{1}{\sin^{2k+1}\nu}\times\nonumber\\ \int_{0}^{\nu}-4k\sin^{2k-1}u\,\cos u\,\int_0^u (\sum_{m=1}^{k-1}c_m(2k-1)\,(\dfrac{\sin v}{\sin u})^{2m})\cos v\,\overline{a_{1}(v)}\, dv\,du=\nonumber\\\overline{a_{1}(u)}+\frac{1}{\sin\nu}\int_0^\nu (\sum_{m=1}^{k-1}c_m(2k-1)(\dfrac{\sin u}{\sin\nu})^{2m}-4k(\dfrac{\sin u}{\sin\nu})^{2k})\cos u\,\overline{a_{1}(u)} du\nonumber\\+\frac{1}{\sin\nu}\int_0^\nu(\sum_{m=1}^{k-1}\dfrac{-4kc_m(2k-1)}{2k-2m+1}[(\dfrac{\sin u}{\sin\nu})^{2m}-(\dfrac{\sin u}{\sin\nu})^{2k}]\cos u\,\overline{a_{1}(u)} du.\nonumber
\end{eqnarray}
Theorem 5 is proved.

\section{The Inversion formula }\label{sec6}

\noindent Using the expressions for $\overline{a_{2k}(\nu)}$ and $\overline{a_{2k+1}(\nu)}$ given in (\ref{6.2}) and (\ref{6.5}) for the partial sum of (\ref{5.2}), we obtain
\begin{eqnarray}
	\overline{F(\pi/2)}+\sum_{k=1}^{n}(\overline{a_{2k-1}(\pi/2)}+\overline{a_{2k}(\pi/2)})=
	\overline{F(\pi/2)}+n(\overline{a_{1}(\pi/2)}+2\overline{a_{0}(\frac{\pi}{2})})+\label{7.1}\\\sum_{k=1}^{n}\int_0^{\pi/2}\left( (\sum_{m=1}^{k-1}c_m(2k-1){\sin^{2m}u})\overline{a_{1}(u)}+ (\sum_{m=1}^{k}c_m(2k)\,{\sin^{2m-1}u})\overline{a_{0}(u)}\right)\cos u\,du\nonumber\\= \overline{F(\pi/2)}+n(\overline{a_{1}(\pi/2)}+2\overline{a_{0}(\frac{\pi}{2})})+\int_0^{\pi/2}\left( P_n^1(u)\overline{a_{1}(u)}+ P_n^0(u)\overline{a_{0}(u)}\right)\cos u\,du
	.\nonumber\end{eqnarray}
Here for $n\geq1$ 
\begin{equation}\label{7.2}
	P_n^1(u)= \sum_{k=1}^{n}\sum_{m=1}^{k-1}c_m(2k-1){\sin^{2m}u},\,\,\,\text{and}\,\,\,P_n^0(u)=\sum_{k=1}^{n}\sum_{m=1}^{k}c_m(2k){\sin^{2m-1}u}\end{equation}
are trigonometric polinomials. 

\noindent Substituting (\ref{7.1}), (\ref{7.2}) into  (\ref{5.2}), we obtain
(\ref{1.10}). Theorem 2 is proved.

\begin{figure}
	\center
	\includegraphics{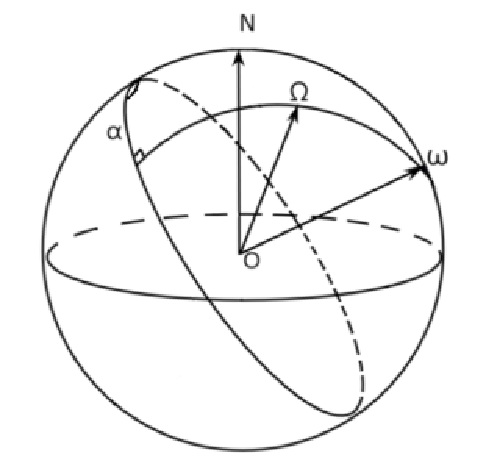}
	\caption{An illustration to Theorem 1.3.}
\end{figure}

\section{The weighted Funk transforms with respect to $\Omega$}\label{sec6}
Let $Cf(\omega)$ and $Sf(\omega)$ be the two weighted Funk transforms of a function $f\in{\textit{C}}^{1}(\mathbf S^2)$ with respect to the North Pole. For $\Omega\in\mathbf S^2$, we have (see Fig.2)
\begin{eqnarray}
	C_\Omega f(\omega)=\dfrac{1}{2\pi}\int_{{\mathbf  S}_\omega}\cos\varphi \,\,f_\omega(\varphi) \,d\varphi=\dfrac{1}{2\pi}\int_{{\mathbf  S}_\omega}\cos(\varphi'-\alpha) \,f_\omega(\varphi) \,d\varphi=\nonumber \\\cos\alpha\, Cf(\omega)+\sin\alpha \,Sf(\omega) 
	\label{7.3}		\end{eqnarray}
here we measure $\varphi$ from ${\Omega}_\omega$ (the projection direction on $\Omega$ to $\mathbf S_\omega$ is the reference point) and correspondingly we measure $\varphi'$ from ${\emph{N}}_\omega$ (the projection direction of $\emph{N}$ to $\mathbf S_\omega$), $\alpha$ is the angular coordinate of ${\Omega}_\omega$ with respect to ${\emph{N}}_\omega$.	Theorem 3 is proved.

\section{Computational Implementation}\label{sec8}

\noindent Now we are going to implement the  iterative reconstruction algorithm obtained from (\ref{1.6}).

\noindent	\textbf{Example}. Consider the following odd function defined on  $\mathbf S^2$

\begin{equation}\label{8.1}
	f(\omega)=	f(\nu,\tau)	=	\cos^3\nu.
\end{equation}
here we use the
standard spherical coordinates $\omega=(\nu,\tau)$ for $\omega\in\mathbf S^2$ determined by choice of the North Pole ${\emph{N}}$.

\noindent	From (\ref{1.3}) it follows (using rules for the right spherical triangle we have $f_(\nu,\tau)(\varphi)=\sin^3\nu\cos^3\varphi$) 
\begin{equation}\label{8.2}
	2\pi	Cf(\omega)=\int_{{\mathbf  S}_(\nu,\tau)}\cos\varphi \,\,f_\omega(\varphi) \,d\varphi=\int_{-\pi}^\pi\cos\varphi \,\,\sin^3\nu\cos^3\varphi\,d\varphi=\frac{3\pi}{4}\sin^3\nu, 
\end{equation}
hence 
\begin{equation}\label{8.2.1}
	\overline{Cf(u)}=\int_0^{2\pi}C f(u,\tau)\,d\tau=\frac{3\pi}{4}\sin^3u, 
\end{equation}

\noindent From (\ref{6.6}) it follows that 
\begin{equation}\label{8.3}
	P_2^1(u)=-4\sin^2 u. 
\end{equation}

\noindent Now we approximate $f$ for $n=2$ using the reconstruction algorithm obtained from (\ref{1.6}).

\noindent Substituting (\ref{8.2}) and (\ref{8.3}) into (\ref{1.10}), and taking into account that $\overline{Ff(u)}=0$,  we obtain the approximation of $f$ for $n=2$

\begin{equation}\label{8.4}
	2\pi=2\pi f(\emph{N})=
	\left(4\overline{Cf(\pi/2)}+\int_0^{\pi/2}\left( -6\pi\sin^2 u\sin^3u\right)\cos u\,du\right)=2\pi.		\end{equation}

\noindent	We obtain the exact reconstruction for n = 2 since $a_{2k-1}(\nu,\tau)=0$ for
any $k>2$  ($\cos^3\varphi$ is orthogonal to $\cos{(2k-1)\varphi}$ for $k>2$).

\section{Conclusion}

It is known that the Funk transform (the Funk-Radon transform)
is invertible in the class of even (symmetric) continuous functions defined on the unit 2-sphere $\mathbf S^{2}$.
In this article, for the reconstruction of 			$f\in C(\mathbf S^{2})$ (can be non-even), an additional condition (to reconstruct an odd function) is found, and the injectivity of the so-called two data Funk transform is considered. Applying the consistency method, suggested by the author of this article, an iterative  inversion formula of the transform is presented. Such inversions have theoretical significance in convexity theory and integral geometry. Also, the Funk-Radon transform is used in Diffusion-weighted magnetic resonance imaging.

\noindent\textbf{ Data availability.} The data that support the findings of this study are available upon research.

\noindent \textbf{Conflict of interest}. The author declares no competing interests.

\end{document}